\documentclass[12 pt]{amsart}

\usepackage{layout}
\usepackage{amsmath, amsthm}
\usepackage{amssymb, amsrefs}
\usepackage{stmaryrd}
\usepackage{graphicx}
\usepackage{setspace}
\usepackage{textcomp}
\usepackage{esvect}
\usepackage[letterpaper]{geometry}
\usepackage{arydshln}
\usepackage[all]{xy}
\usepackage{etex}
\usepackage{ifthen}

\usepackage[dvipsnames]{xcolor}

\usepackage{epsfig}
\usepackage{latexsym}
\usepackage{amsthm}

\usepackage{mathrsfs}

\usepackage[colorlinks,citecolor=blue]{hyperref}

\usepackage[latin1]{inputenc}

\usepackage{tikz-cd}
\usetikzlibrary{graphs,graphs.standard}
\usepackage{pgfplots}
\usetikzlibrary{matrix,arrows,decorations.pathmorphing}
\usetikzlibrary{shapes.geometric}

\newtheorem{theorem}{Theorem}[section]
\newtheorem{lemma}[theorem]{Lemma}
\newtheorem{proposition}[theorem]{Proposition}

\newtheorem{corollary}[theorem]{Corollary}

\theoremstyle{definition}
\newtheorem{definition}[theorem]{Definition}

\newtheorem{example}[theorem]{Example}

\newtheorem{obs}[theorem]{Observation}
\newtheorem{notation}[theorem]{Notation}

\theoremstyle{definition}
\numberwithin{equation}{section}

\tikzset{
    ncbar angle/.initial=90,
    ncbar/.style={
        to path=(\tikztostart)
        -- ($(\tikztostart)!#1!\pgfkeysvalueof{/tikz/ncbar angle}:(\tikztotarget)$)
        -- ($(\tikztotarget)!($(\tikztostart)!#1!\pgfkeysvalueof{/tikz/ncbar angle}:(\tikztotarget)$)!\pgfkeysvalueof{/tikz/ncbar angle}:(\tikztostart)$)
        -- (\tikztotarget)
    },
    ncbar/.default=0.5cm,
}

\tikzset{square left brace/.style={ncbar=0.5cm}}
\tikzset{square right brace/.style={ncbar=-0.5cm}}

\tikzset{round left paren/.style={ncbar=0.5cm,out=120,in=-120}}
\tikzset{round right paren/.style={ncbar=0.5cm,out=60,in=-60}}

\DeclareMathOperator{\Aut}{Aut}

\DeclareMathOperator{\zed}{\mathbb{Z}}

\def \W {\mathfrak{W}}

\def \Gph {\mathsf{Gph}}

\def \W {\mathfrak{W}}

\definecolor{laura}{rgb}{.4, 0, .6}

\begin{document}

 \title{Homotopy Covers of Graphs}
    \author{T. Chih and L. Scull}

\maketitle

\vspace*{-2em}



\begin{abstract}
We develop a theory of $\times$-homotopy, fundamental groupoids and covering spaces that applies to non-simple graphs, generalizing existing results for simple graphs.   We prove that $\times$-homotopies from finite graphs can be decomposed into moves that adjust at most one vertex at a time, generalizing the spider lemma of \cite{CS1}.   We define a notion of homotopy covering map and develop a theory of universal covers and deck transformations, generalizing  \cite{Matsushita}, \cite{TardifWroncha} to non-simple graphs.  We examine the case of reflexive graphs (each vertex having at least one loop).     We also  prove that these homotopy covering maps satisfy a homotopy lifting property for arbitrary graph homomorphisms,  generalizing path lifting results of \cite{Matsushita}, \cite{TardifWroncha}.

\end{abstract}

\section{Introduction}

 Covers of graphs were originally studied by viewing graphs as   1-dimensional topological spaces.  From this perspective, many properties of basic topology have been extended to graphs, including the fundamental groupoid, covering spaces, universal covers and deck transformations  \cite{Angluin1}, \cite{Bass},  \cite{KwakNedela}.   These results have been developed for non-simple graphs that allow loops and parallel edges.  
 
Viewed as 1-dimensional topological spaces, the homotopies between graphs are limited.   However, it is possible to develop theories of homotopies between graph homomorphisms that allow more interesting deformations and take into account more of the graph structure.   
   There are two prominent theories of homotopies for graphs in the literature:     $A$-homotopy \cite{Babson1}, \cite{Barcelo1}, \cite{BoxHomotopy},  \cite{hardeman2019lifting} and $\times$-homotopy \cite{CompGH}, \cite{CS1}, \cite{CS2}, \cite{Docht1}, \cite{Docht2}, \cite{Kozlov1}, \cite{Kozlov3}, \cite{KozlovShort}.   Covers have been considered in the context of  $A$-homotopy in \cite{hardeman2019lifting}, where  Hardeman  shows that for graphs with no three or four cycles, graph covers lift  $A$-homotopies.

   The goal of this paper is to  understand what $\times$-homotopy, the fundamental group(oid) and covering maps of graphs look like in the context of graphs which allow loops and multiple parallel edges.     We begin with some basic results about homotopies of graph homomorphisms between graphs with parallel edges in Section \ref{S:homotopy}, showing that the spider lemma of \cite{CS1} still holds, allowing us to decompose a $\times$-homotopy between finite graphs to a sequence of spider moves that adjust at most  one vertex at a time.     We then use this to define the fundamental groupoid of non-simple graphs in Section \ref{S:groupoid}, following \cite{CS2}, and give a relationship between the fundamental groupoid of a reflexive graph and its unlooped counterpart in Propositon  \ref{P:reflexivegroupoid}.      In Section \ref{S:HC}, we define homotopy covering maps and develop their properties, including a theory of deck transformations, for non-simple graphs.  We also examine in detail the case of reflexive graphs.    We conclude this section with proving that  these homotopy covers satisfy the general homotopy lifting property for any graph homomorphism in  Theorem \ref{T:HLiftGeneral}.  This last result generalizes results in the literature in two ways: in addition to extending to non-simple graphs, our result allows lifting of homotopies between any graph homomorphisms, and not just paths as in \cite{Matsushita}, \cite{TardifWroncha}.

   When restricted to simple graphs, our results recover several known results in the literature.         Our fundamental groupoid recovers that of Chih-Scull \cite{CS2} when we allow loops but not parallel edges, and that of  Sankar \cite{Sankar} for simple graphs.    Additionally, the spider lemma shows that the definition we use of $\times$-homotopy is equivalent to the definition of $r$-homotopy in Matsushita \cite{Matsushita} when $r=2$, as this definition allows homomorphisms to be shifted by $r$ edges (and hence $r-1$ vertices).   For simple graphs,    Matsushita develops covering space theory for $r$-homotopy based on a neighborhood complex.   Lastly, for simple graphs Tardif-Wroncha \cite{TardifWroncha} have a theory of $\diamond$-covers which lift any embedded 4-cycle, and in the appendix they give an alternate formulation which matches up with our Definition \ref{D:hocover} of homotopy covers based on 2-neighborhoods.  
   
\section{Graphs and Homotopy} \label{S:homotopy}

We begin with the definition of our category of graphs, following \cite{Bass}, \cite{KwakNedela}.

\begin{definition} The category of graphs $\Gph$ is defined by: \begin{itemize} \item An object is a  graph $G$,  consisting of a  set of vertices $V(G) = \{ v_{\lambda}\}$ and a set $E(G) = \{ e_\beta \} $ of edges, with maps $\partial_0, \partial_1:  E(G) \to V(G)$ and an involution $e \to \overline{e}$ of $E(G)$ such that $\partial_i \overline{e}= \partial_{1-i} e$.       
\item 
A homomorphism of graphs $f:  G \to H$ is given by a vertex map $f:  V(G) \to V(H)$ and an edge map $f:  E(G) \to E(H)$ such that $f(\partial_i e) =\partial_i f(e)$ and $f(\overline{e}) = \overline{f(e)}$. 
 \end{itemize} \end{definition}  

   This category of graphs allows loops and  multiple parallel edges between the same vertices. Bass \cite{Bass} requires that the involution be fixed point free, while Kwak-Nedela \cite{KwakNedela}  allows 'semi-edges' such that $\overline{e} = e$. 

In considering graphs which allow both loops and parallel edges, we compare these to their single-edge and unlooped versions.   Thus we will use the following notation.

\begin{notation} \label{N:unloop}
If $G$ is a graph, then:
\begin{itemize}
    \item $G_u$ refers to the graph $G$ with all loops deleted:  $V(G_u) = V(G)$ and \\  $E(G_u)= E(G) \backslash \{e : \partial_0 (d) = \partial_1(d)\}$.
    \item $G_s$ refers to the single-edge graph with parallel edges have been identified:  $V(G_s) = V(G)$ and $E(G_s) =\{ [e] : e \in E(G) \textup{ and } [e] = [e']  \textup{ if } \partial_0(e) = \partial_0(e') \textup{ and } \partial_1(e) = \partial_1(e') \} $.
\end{itemize}
    
\end{notation}

In this paper, we consider homotopy between  graph homomorphisms described above, defined using the categorical product.   This is  sometimes referred to as $\times$-homotopy.   This is the only version of homotopy that will appear in this paper, and we will  refer to it simply as `homotopy'.   It uses the following definition.  

\begin{definition}
    The categorical product $G \times H$ is defined by  $V(G \times H) = V(G) \times V(H)$ and $E(G \times H) = E(G) \times E(H)$  with $\partial_i(e_1, e_2) = (\partial_i(e_1), \partial_i(e_2))$ and $\overline{(e_1, e_2)} = (\overline{e}_1, \overline{e}_2)$.
   
\end{definition}

We define homotopy using the product with the interval graph $I_n$.  
\begin{definition} \label{D:int}
     Let $I_n$ denote the looped path graph with vertices $V(I_n) = \{ 0, 1, 2, \dots, n\}$ and edges $E(I_n) = \{ \ell_i, e_i, \overline{e}_i\}$ for $1 \leq i \leq n$, where $\partial_i(\ell_j ) = j$ and $\partial_0(e_j) = j-1, \partial_1(e_j) = j$. 
\end{definition}

Observe that for any $j$, the vertices of the form $(v, j)$ and edges of the form $(e, \ell_j)$ define a subgraph of $G \times I_n$ which is isomorphic to $G$;  we denote this subgraph $G \times \{j\}$.     The following definition is modeled on that of \cite{Docht1}.

\begin{definition}  \label{D:htpy}  For $f, g:  G \to H$, we say that $f$ is {\bf homotopic} to $g$, written $f \simeq g$,  if there is a graph homomorphism  $\Lambda: G \times I_n  \to H$ such that $\Lambda | _{G \times \{ 0\} } = f$ and $\Lambda | _{G \times \{ n\} } = g$.   We write $f \simeq g$.         \end{definition}

The introduction of multiple edges allows for multiple graph homomorphisms which agree on vertex sets.   We show that these are all homotopic.

\begin{lemma} \label{L:hom}
    If $f, g:  G \to H$ are graph homomorphism which agree on vertices, then  $f \simeq g$.  
\end{lemma}

\begin{proof}
    We can define a homotopy $\Lambda:  G \times I_1 \to H$ by $\Lambda( x, i) = f(x) = g(x) $ and 
       \begin{equation*}
{\Lambda}(e, e')= \begin{cases} f(e)  & \textup{ if } e' = \ell_0\\ g(e)  & \textup{ otherwise } 
\end{cases}  
\end{equation*} 
     
\end{proof}

We can use this to show that for a loop-free graph $G$,  the graph $G$ is homotopy equivalent to the graph $G_s$   where all parallel edges have been identified. 
\begin{proposition}
    Let $G$ be a loop-free graph, and $G_s$ the single-edge graph as in Notation \ref{N:unloop}.   Then $G_s$ is a strong deformation retract of $G$. 
\end{proposition}

\begin{proof}
    Define the projection $\pi:  G \to G_s$ by the identity on vertices and the projection on edges.  Define $\iota:  G_s \to G$ using the identity on vertices, and  for each pair of edges $[e], [\overline{e}]$ choose an edge $e$ between the same vertices and define $\iota[e] = e$ and $\iota[\overline{e}] = \overline{e}$.    Then $\iota \pi = id_H$ and $\pi \iota \simeq id_G$ by Lemma \ref{L:hom}.
\end{proof}

Thus in any homotopy invariant construction of a loop-free graph, we may substitute $G_s$ instead of $G$.

When working with homotopies between graph homomorphisms  that do not necessarily agree on vertices, the following is a useful tool.   

\begin{proposition} \label{P:spider} If $G$ is finite and  $f, g:  G \to H$ then  $f \simeq g$ if and only if there is a finite sequence of spider moves connecting $f$ and $g$, where each spider move changes the value of $f$ by at most a  single vertex, with the additional condition that if vertex $v$ has a loop, then any spider move that changes the value on $v$ must change to a connected vertex. 
\end{proposition}

\begin{proof} The proof of \cite{CS1} Proposition 4.4 can by adapted to the multiple-edge case as follows. If $f$ and $g$ agree on vertices,  then by Lemma \ref{L:hom} there is a one step homotopy between them.   If $f, g$ differ by a single vertex $v$ then we can define a homotopy $\Lambda:  G \times I_1 \to H$ by 
$\Lambda | _{G \times \{ 0\} } = f$ and $\Lambda | _{G \times \{ 1\} } = g$, and for a non-loop edge $e$, we define
\begin{equation*}
{\Lambda}(e, e_1) = \begin{cases} f(e)  & \textup{ if } \partial_0e = v \\ g(e)  & \textup{ if } \partial_1e = v
\end{cases}
\end{equation*} and \begin{equation*}
{\Lambda}(e, \overline{e}_1) = \begin{cases} f(e)  & \textup{ if } \partial_1e = v \\ g(e)  & \textup{ if } \partial_0e=v
\end{cases}
\end{equation*}
If $v$ has any loops $\ell_i$, then  we know there is an edge $\eta$ connecting $f(v)$ to $g(v)$ and we define $\Lambda (\ell_i, e_1) = \eta$ and $\Lambda (\ell_i, \overline{e}_1) = \overline{\eta}$.

Conversely, suppose we have a  homotopy $f \simeq g$; it is sufficient to consider a one-step homotopy   $\Lambda:  G \times I_1 \to H$ between $f$ and $g$ and then induct to get longer homotopies.    Suppose that  $G$ has vertex set $\{ v_1, v_2, \dots, v_n\}$.   We define a homotopy $\widetilde{\Lambda}: G \times I_n \to H$ on vertices by 
\begin{equation*}
\widetilde{\Lambda}(v_i, k) = \begin{cases} \Lambda(v_i, 0) = f(v_i) & \textup{ if } i>k \\ \Lambda(v_i, 1)= g(v_i) & \textup{ if } i\leq k
\end{cases}
\end{equation*}
 and on edges by the following:  if $e$ is an edge of $G$ with $\partial_0e = v_i$ and $\partial_1e = v_j$ then 
 \begin{equation*}
\widetilde{\Lambda}(e, e') = \begin{cases} \Lambda(e, \ell_0) & \textup{ if } i > {k}, j > {k} \\  \Lambda(e, \ell_1) & \textup{ if } i \leq {k}, j \leq {k} \\  \Lambda(e, e_{1}) & \textup{ if } i > {k}, j \leq {k} \\  \Lambda(e, \overline{e}_{1}) & \textup{ if } i \leq {k}, j > {k} \\ 
\end{cases}
\end{equation*}
 This defines intermediate graph homomorphisms $f_k = \widetilde{\Lambda}(v, k)$ in which $f_k, f_{k+1}$ either agree on vertices or differ by a single vertex, and thus we get a sequence of spider moves joining $f$ to $g$.       \end{proof}

\section{The Fundamental Groupoid}\label{S:groupoid}
In this section, we review the fundamental groupoid of \cite{CS2} and extend it to the multiple-edge case.  We then consider the particular case of reflexive graphs. 

\begin{definition}\label{d:walk}
    A walk of length $n$ is defined by a sequence of edges $(e_1e_2 \dots e_n)$ such that $\partial_1(e_i) = \partial_0(e_{i+1})$.   If $\partial_0(e_1) = v$ and $\partial_1(e_n) =w$ we say that the walk starts at $v$ and ends at $w$.  
\end{definition}

This is used to define what we refer to as the walk groupoid (called the fundamental groupoid in \cite{Bass}, \cite{KwakNedela}).  We will present this in terms of the concept of `prunes' from \cite{CS2}.  
\begin{definition} 
    When a walk $\alpha = (e_1 e_2 \dots e_n)$ satisfies $e_{i-1} = \overline{e}_i$, we refer to the walk $\alpha'=(e_1 e_2 \dots {e_{i-2} e_{i+1}} \dots e_n)$ that omits the pair ${e_{i-1} e_{i}}$ as a prune of $\alpha$, and $\alpha$ as an anti-prune of $\alpha'$. 
\end{definition}
\begin{definition}
    The walk groupoid $\W(G)$ is a groupoid with: 
    \begin{itemize}
        \item objects given by vertices $v \in V(G)$,
        \item arrows starting $v$ and ending at $w$ given by prune classes of walks, where we identify any walks that are connected by a sequence of prunes and anti-prunes
        \item composition defined by concatenation of walks. 
    \end{itemize}
\end{definition}

This is presented in \cite{Bass} as reduced walks, in which no prunes are possible.  
This walk groupoid does not identify homotopic walks.  
There are several equivalent definitions of a fundamental groupoid which does.   We extend Definition 4.1 of \cite{CS2} to  the following.   

\begin{definition} \label{D:fg}   Let $\Pi(G)$ be the fundamental groupoid of $G$  defined by the following:  

\begin{itemize}
    \item an object of $\Pi(G)$ is a vertex of the graph $G$,
    \item an arrow starting at $v_0$ and ending at $v_n$ in $\Pi(G)$ is given by a prune class of walks with the same start and end, up to  homotopy rel endpoints,
    \item composition of arrows is defined using concatenation of walks. 
\end{itemize}
\end{definition}

The results of Section \ref{S:homotopy} allow us to describe any homotopy of paths as a sequence of spider moves.   Thus our fundamental groupoid consists of equivalence classes of walks, where two walks are equivalent if a sequence of prunes, anti-prunes, and spider moves can transform one into the other.

\begin{notation}\label{N:equiv}
    An arrow in $\Pi(G)$ is defined by a sequence of edges  $(e_1e_2 \dots e_n)$.   However, Lemma \ref{L:hom} ensures that any choice of edges connecting the same sequence of vertices results in the same element of $\Pi(G)$.   Thus in the rest of this section, we will denote arrows in the fundamental groupoid via their sequence of vertices $(v_0 v_1 \dots v_n)$ and not specify edges.  
\end{notation}

We have seen that the addition of parallel  edges does not change the homotopy behavior, so will not be interesting when examining $\Pi(G)$.   However, we are interested in the effect of including loops in our graphs.

\begin{obs}\label{L:loops}  In the fundamental groupoid $\Pi(G)$
\begin{itemize} 
    \item  Any loop $\ell = (vv) $ satisfies $\ell^2 = (vvv) = 1_v$ via a prune.   \item  There is a spider move  from  $(vvw) $ to  $(vww)$.   
    \item If all vertices in a path starting at $v$, and ending at $w$ are looped, we can always move the loops through the path and   write the path as $\gamma * (ww)^k$ for $\gamma$ containing no loops and $k=0 $ or $1$.   
\end{itemize}
    
\end{obs}

In the remainder of this section, we examine the case of reflexive graphs, which we take here to mean that every vertex has at least one loop.  Triangles play a particular role here, so we make the following definition.

\begin{definition}
For a graph $G$ and its fundamental groupoid $\Pi(G)$, we define $\Pi(G)/T$ to be the quotient of the fundamental groupoid $\Pi(G)$ by the relation generated by triangles defined by an embedded $C_3$.   Explicitly, we declare that if we have any embedded $C_3$ defined by $xyzx$, then  if $\alpha = \alpha_1 * \alpha_2$ in $\Pi(G)$ then $\alpha = \alpha_1 * (xyzx) * \alpha_2$ in $\Pi(G)/T$.   Thus we can insert or delete triangles from any path in $\Pi(G)/T$.   
    
\end{definition}
Recall from Notation \ref{N:unloop} that $G_u$ refers to the unlooped graph obtained by removing all loops from $G$.
\begin{proposition}\label{P:reflexivegroupoid}
    Let $G$ be a reflexive graph.   Then $\Pi(G)\cong \Pi(G_u)/T\times \mathbb{Z}/2$.
\end{proposition}

\begin{proof}
The product groupoid $\Pi(G_u)/T \times \zed/2$ is defined as the product of categories, interpreting $\zed/2$ as a category with a single object.  Explicitly, 
 $\Pi(G)/T\times \zed/2$ has the same objects as $\Pi(G)/T$, with arrows $x \to y$ defined by $(\alpha, i)$ for $\alpha : x \to y$ in $\Pi(G)/T$ and $i = 0, 1 \in \zed/2$.   The composition is defined by $(\alpha, i) \circ (\beta, j) = (\alpha \circ \beta, i+j \textup{ (mod 2)} ) $. 
 
   We define  the functor $\Psi:\Pi(G)\to \Pi(G_u)/T\times \mathbb{Z}/2$ by $\Psi(v) = v$ on objects and $\Psi(\alpha) = (\alpha_u, i)$ on arrows, where $\alpha_u$ denotes the walk $\alpha$ with all loops removed, and $i$ is the length of the walk $\alpha$ (mod 2).

    To show that $\Psi$ is well defined, we suppose that $[\alpha] = [\alpha']$ in $\Pi(G)$.  Then $\alpha$ and $\alpha'$ are connected by a sequence of prunes, anti-prunes, and spider moves, all of which preserve the parity of the length.   It is easy to see that any prune, anti-prune, or spider move not involving any loops has a corresponding move on the unlooped version, and that a prune involving a loop must just insert or delete a pair of loops, leaving the unlooped walk the same.   So we consider the case when we have a spider move involving a loop.  Unlooping $(vvw)$ and $(vww)$ both result in $(vw)$.     If we replace $\alpha = uuv$ with $\alpha' = uxv$, then $\alpha_u = uv$ while $\alpha'_u = uxv = uxvuv = uv$ since $uxvu$ is a triangle.   Hence the $\alpha_u = \alpha'_u $ in $\Pi(G_u) / T$.    $\Psi$ respects composition since $(\alpha* \beta)_u = \alpha_u * \beta_u$.

We define an inverse functor $\Phi:\Pi(G_u)/T\times \mathbb{Z}/2 \to  \Pi(G) $ again using the identity on objects, and on arrows defining $\Phi(\gamma, i) = \gamma * (ww)^k$ where $(ww)$ denotes the loop on the ending vertex $w$ of $\gamma$, and $k = i + p $ where  $p$ is the mod 2 length of $\gamma$.  This  ensures that the parity of $\gamma*(ww)^k$ matches the parity of $i$.  

We check that this is well-defined: again, any prune, anti-prune, or spider move of $\gamma$ gives a corresponding move on $\gamma * (ww)$, so we just need to check what happens when we insert or remove a triangle. 
Suppose that we have two walks in $G_u$ defined by $\gamma$ and $\gamma' = \gamma_1 * (xyzx)*\gamma_2$ where $\gamma_1* \gamma_2 = \gamma$, so that we obtain $\gamma'$ by inserting a triangle into $\gamma$.   Then the length of $\gamma$ and $\gamma'$ have opposite parities, and so when we apply $\Phi$ to $(\gamma, i)$ and to $(\gamma', i)$,  one will be concatenated with a loop and the other will not.   Assume that $\Phi(\gamma, i) = \gamma* (ww)$ and $\Phi(\gamma', i) = \gamma'$.   Then we use Observation \ref{L:loops}  to get the following in $\Pi(G)$:   $\gamma'*(ww) = \gamma_1 * (xyzx)*\gamma_2* (ww) = \gamma_1 * (xyzx)*(xx) *\gamma_2 = \gamma_1 *(xyzxx) * \gamma_2 = \gamma_1 * (xyxxx) * \gamma_2 = \gamma_1 * (xyx) * \gamma_2 = \gamma_1 *\gamma_2 = \gamma$.   Therefore $\gamma'  = \gamma'* (ww) * (ww)  = \gamma*(ww)$.   A similar argument verifies well-definedness when $\Phi$ puts the loop on $\gamma'$, and migrating loops through out paths ensures that $\Phi$ respects the concatenation operation.   

 Because all loops can be moved to the end of any walk in $\Pi(G)$, $\Psi $ and $\Phi $ are inverse functors, giving the desired isomorphism of groupoids.

\end{proof}

\begin{corollary}\label{C:reflexivePi}
 Suppose that $G$ is a reflexive graph with no embbedded $C_3$.   Then the fundamental groupoid  $\Pi(G) \simeq \Pi(G_u) \times \zed/2$.
\end{corollary}

\begin{example}\label{E:reflexivegroupoid}
    Consider the following graph $G$ and its  unlooped subgraph $G_u$: 

    \[\begin{tikzpicture}
\draw ({sin(36+144)},{cos(36+144)})\foreach \x in {3,...,6}{--({sin(72*\x+36)}, {cos(72*\x+36)})};

\draw ({sin(36)},{cos(36)})--(0,1.51) -- ({-1*sin(36)},{cos(36)});

\draw ({sin(36+72)},{cos(36+72)}) edge[bend left] ({sin(36+144)},{cos(36+144)});
\draw ({sin(36+72)},{cos(36+72)}) edge[bend right] ({sin(36+144)},{cos(36+144)});

\draw[fill, ] ({sin(36+72*0)}, {cos(36+72*0)}) node[right]{$a$}  circle (2pt);
\draw[fill, ] ({sin(36+72*1)}, {cos(36+72*1)}) node[below right]{$b$}  circle (2pt);
\draw[fill,   ] ({sin(36+72*2)}, {cos(36+72*2)}) node[below]{$c$}  circle (2pt);
\draw[fill, ] ({sin(36+72*3)}, {cos(36+72*3)}) node[below left]{$d$}  circle (2pt);
\draw[fill,  ] ({sin(36+72*4)}, {cos(36+72*4)}) node[left]{$e$}  circle (2pt);
\draw[fill] (0, 1.51)  node[right]{$x$}  circle (2pt);

\draw (0.588, 0.809)  to[in=50,out=140,loop] (0.588, 0.809);
\draw (0.951, -0.309)  to[in=50,out=140,loop] (0.951, -0.309);
\draw (0,-1)  to[in=50,out=140,loop] (0,-1);
\draw (-0.588, 0.809)  to[in=50,out=140,loop] (-0.588, 0.809);
\draw (-0.951, -0.309)  to[in=50,out=140,loop] (-0.951, -0.309);
\draw (0,1.51)  to[in=50,out=140,loop] (0,1.51);

\node at (0, -1.5){$G$};

\end{tikzpicture}\ \ \ \ \ \ \ \ \ \ \ \ \ \ \ \ \ \ \ \ \ 
\begin{tikzpicture}
\draw ({sin(36+144)},{cos(36+144)})\foreach \x in {3,...,6}{--({sin(72*\x+36)}, {cos(72*\x+36)})};

\draw ({sin(36)},{cos(36)})--(0,1.51) -- ({-1*sin(36)},{cos(36)});

\draw ({sin(36+72)},{cos(36+72)}) edge[bend left] ({sin(36+144)},{cos(36+144)});
\draw ({sin(36+72)},{cos(36+72)}) edge[bend right] ({sin(36+144)},{cos(36+144)});

\draw[fill, ] ({sin(36+72*0)}, {cos(36+72*0)}) node[right]{$a$}  circle (2pt);
\draw[fill, ] ({sin(36+72*1)}, {cos(36+72*1)}) node[below right]{$b$}  circle (2pt);
\draw[fill,   ] ({sin(36+72*2)}, {cos(36+72*2)}) node[below]{$c$}  circle (2pt);
\draw[fill, ] ({sin(36+72*3)}, {cos(36+72*3)}) node[below left]{$d$}  circle (2pt);
\draw[fill,  ] ({sin(36+72*4)}, {cos(36+72*4)}) node[left]{$e$}  circle (2pt);
\draw[fill] (0, 1.51)  node[right]{$x$}  circle (2pt);

\node at (0, -1.5){$G_{u}$};

\end{tikzpicture}
\]

  The arrows of $\Pi G_u$ are walks in $G_u$ up to edge shuffles, since these  are the only possible spider moves.  Therefore the closed walks  $(aexa), (abcdea)$ do not commute in $\Pi G_u$.   When we add the loops and move to the graph $G$, we introduce loops into our walks and the relations  $(aexa)=(aa), (vv)^2=(v)$ and  $(uuv)=(uvv)$.  Thus we have commuting closed walks $(abcde)$ and  $(aexa)=(aa)$  and  $\Pi G\cong \Pi G_u/\langle (aexa) \rangle \times \mathbb{Z}/2$.
\end{example}

\section{Homotopy Covers of Graphs}\label{S:HC}

This section builds on the fundamental groupoids of Section \ref{S:groupoid} to develop the theory of homotopy covers of graphs.   We define homotopy covering maps and  show that they satisfy a version of the classical theory of universal covers and deck transformations, and describe the homotopy covers of reflexive graphs.   We close by showing that these maps satisfy the homotopy lifting property for any graph homomorphisms in Theorem \ref{T:HLiftGeneral}.  Throughout this section, we will assume that all graphs are connected. 

We begin with the classical definition of graph covers. 
  \begin{definition}\label{D:Cover}  \cite{KwakNedela}
A covering map is a  graph morphism $f:\widetilde{G}\to G$ such that
\begin{itemize} 
\item $f$ is a surjection on vertices and edges,

\item  $f$ induces a bijection on neighborhoods:   let $N_E(v)$ denote the set of edges $e$ with $\partial_0(e)=v$.  Then for any vertex $\tilde{v}$, with $f(\tilde{v}) = v$,   $f$ induces a bijection  $N_E(\tilde{v})\to N_E(f(\tilde{v}))$.

\end{itemize}

\end{definition}

This coincides with the definition given in \cite{Angluin1} for simple graphs,  since in that case  for $e_1, e_2\in N_E(v), \partial_1(e_1)=\partial_1(e_2)$ if and only if $e_1=e_2$.

An important property of covering maps is the following path lifting result.  
\begin{lemma}\cite{Angluin1}, \cite{Hatcher}, \cite{KwakNedela}\label{L:liftwalkold}
Suppose that  $f:\widetilde{G}\to G$ is a covering map.  Then given  a walk $\alpha$ in $G$  starting at $v$ defined by  $\alpha=(e_1e_2\cdots e_n)$    and any $\tilde{v}\in f^{-1}(v)$, there exists a unique walk $\tilde{\alpha}$  in $\widetilde{G}$ starting at $\tilde{v}$ defined by  $\tilde{\alpha}=(d_1d_2\cdots d_n)$ such that $f(\tilde{\alpha})=\alpha$.
\end{lemma}
 
The covers of Definition \ref{D:Cover} do not incorporate the notion of homotopy of graphs.  Thus, we consider the following adaptation, which requires bijections on 2-neighborhoods rather than only neighborhoods.   

\begin{definition}\label{D:2Nbd}
For any vertex $v \in V(G)$ we define the extended neighborhood $N_2(v)$ to be the walks of length $2$ which start at $v$.  
\end{definition}

Our homotopy covering  maps induce bijections on these 2-neighborhoods.  

\begin{definition}\label{D:hocover}

A surjective graph morphism between  graphs  $f:\widetilde{G}\to G$ is a \textit{homotopy covering map} if given any $v\in V(G)$ and $\tilde{v}\in V(\widetilde{G})$ such that $\tilde{f}(\tilde{v})=v$, then $f$ induces a bijection $N_2(\tilde{v}) \to N_2 (v)$.  Furthermore, this bijection respects endpoints in the sense that two walks in $N_2(\tilde{v}), (d_1d_2), (d'_1d'_2)$ have $\partial_1(d_2)=\partial_1(d'_2)$ in $\widetilde{G}$ if and only if $\partial_1(f(d_2))=\partial_1(f(d'_2))$ in $G$.  \end{definition}

In working with these homotopy covering maps, the following elementary observations will be extremely useful. 

\begin{obs} \label{O:useful}

    Suppose that $f:  \widetilde{G} \to G$ is a homotopy covering map. 
    \begin{itemize}
        \item If a 2-walk in $N_2(v)$ is of the form  $e \overline{e}$, then the lift must have the same form, since if we have a lift $d d'$  then  $d \overline{d}$  defines a 2-walk which also projects down to $e \overline{e}$ and so by uniqueness of lifting we must have $d' = \overline{d}$.  
        \item Any homotopy cover is a cover, since $N_E(v)$ is in bijection with  2-walks of the form $e \overline{e}$. 
        \item Parallel edges lift to parallel edges, since if $e_1, e_2$  both have $\partial_0(e_i) = v$ and $\partial_1(e_i) = w$ then $e_1 \overline{e}_2$ is a 2-walk and so lifts to a 2-walk $d d'$ from $\widetilde{v}$ to $\widetilde{v}$.   Then $\overline{d'}$ is the unique lift of $e_2$.  
    \end{itemize}
\end{obs}

\begin{example}\label{Example:HCover} 
The following is a homotopy covering map taking $v_i \to v$. 
Here we see that 2-walks  $abc, ab'c$ with the same endpoint lift to 2-walks with the same endpoint  in $\widetilde{G}$.

$$
\begin{tikzpicture}[scale=.7]
\node(G) at (0,0){\begin{tikzpicture}[scale=.7]
\draw ({sin(144)},{cos(144)})\foreach \x in {3,...,6}{--({sin(72*\x)}, {cos(72*\x)})};
\draw ({sin(144)},{cos(144)}) edge[bend left] ({sin(72)},{cos(72)});
\draw ({sin(144)},{cos(144)}) edge[bend right] ({sin(72)},{cos(72)});
\draw (-.951, 0.309)--(0,0)--(.951,0.309);

\draw[ ] ({sin(72*0)}, {cos(72*0)}) -- node[above]{$b$} ({sin(72*0)}, {cos(72*0)});
\draw[  ] ({sin(72*1)}, {cos(72*1)}) -- node[right]{$c$} ({sin(72*1)}, {cos(72*1)});
\draw[ ] ({sin(72*2)}, {cos(72*2)}) -- node[below right]{$d$} ({sin(72*2)}, {cos(72*2)});
\draw[ ] ({sin(72*3)}, {cos(72*3)}) -- node[below left]{$e$} ({sin(72*3)}, {cos(72*3)});
\draw[ ] ({sin(72*4)}, {cos(72*4)}) -- node[left]{$a$} ({sin(72*4)}, {cos(72*4)});
\draw[  ] (0, 0) -- node[below]{$b'$} (0, 0);

\draw[fill,  ] ({sin(72*0)}, {cos(72*0)}) circle (2pt);
\draw[fill,   ] ({sin(72*1)}, {cos(72*1)}) circle (2pt);
\draw[fill,  ] ({sin(72*2)}, {cos(72*2)}) circle (2pt);
\draw[fill,  ] ({sin(72*3)}, {cos(72*3)}) circle (2pt);
\draw[fill,  ] ({sin(72*4)}, {cos(72*4)}) circle (2pt);
\draw[fill,   ] (0, 0) circle (2pt);

\node at (0, -1.5){$G$};

\end{tikzpicture}};

\node(T) at (-8,0){\begin{tikzpicture}[scale=0.7]
\draw ({2*sin(36*2)},{2*cos(36*2)})\foreach \x in {3,...,6}{--({2*sin(36*\x)}, {2*cos(36*\x)})};
\draw ({2*sin(7*36)},{2*cos(7*36)})\foreach \x in {8,...,11}{--({2*sin(36*\x)}, {2*cos(36*\x)})};
\draw ({2*sin(36)},{2*cos(36)}) edge[bend left]  ({2*sin(36*2)},{2*cos(36*2)});
\draw ({2*sin(36)},{2*cos(36)}) edge[bend right]  ({2*sin(36*2)},{2*cos(36*2)});
\draw ({2*sin(36*6)},{2*cos(36*6)}) edge[bend left]  ({2*sin(36*7)},{2*cos(36*7)});
\draw ({2*sin(36*6)},{2*cos(36*6)}) edge[bend right]  ({2*sin(36*7)},{2*cos(36*7)});
\draw (-1.176, 1.618)--(0,1.236)--(1.176,1.618);
\draw (-1.176, -1.618)--(0,-1.236)--(1.176,-1.618);

\draw[ ] ({2*sin(36*0)}, {2*cos(36*0)}) -- node[above]{$b_1$} ({2*sin(36*0)}, {2*cos(36*0)});
\draw[  ] ({2*sin(36*1)}, {2*cos(36*1)}) -- node[above right]{$c_1$} ({2*sin(36*1)}, {2*cos(36*1)});
\draw[ ] ({2*sin(36*2)}, {2*cos(36*2)}) -- node[right]{$d_1$} ({2*sin(36*2)}, {2*cos(36*2)});
\draw[ ] ({2*sin(36*3)}, {2*cos(36*3)}) -- node[right]{$e_1$} ({2*sin(36*3)}, {2*cos(36*3)});
\draw[ ] ({2*sin(36*4)}, {2*cos(36*4)}) -- node[below right]{$a_2$} ({2*sin(36*4)}, {2*cos(36*4)});
\draw[ ] ({2*sin(36*5)}, {2*cos(36*5)}) -- node[below]{$b_2$} ({2*sin(36*5)}, {2*cos(36*5)});
\draw[  ] ({2*sin(36*6)}, {2*cos(36*6)}) -- node[below left]{$c_2$} ({2*sin(36*6)}, {2*cos(36*6)});
\draw[ ] ({2*sin(36*7)}, {2*cos(36*7)}) -- node[left]{$d_2$} ({2*sin(36*7)}, {2*cos(36*7)});
\draw[ ] ({2*sin(36*8)}, {2*cos(36*8)}) -- node[left]{$e_2$} ({2*sin(36*8)}, {2*cos(36*8)});
\draw[ ] ({2*sin(36*9)}, {2*cos(36*9)}) -- node[above left]{$a_1$} ({2*sin(36*9)}, {2*cos(36*9)});
\draw[  ] (0,1.236) -- node[below]{$b'_1$}(0,1.236);
\draw[  ] (0,-1.236) -- node[above]{$b'_2$}(0,-1.236);

\foreach \x in {0,...,1}{\draw[fill,  ] ({2*sin(180*\x)}, {2*cos(180*\x)}) circle (2pt);}
\foreach \x in {0,...,1}{\draw[fill,   ] ({2*sin(180*\x+36)}, {2*cos(180*\x+36)}) circle (2pt);}
\foreach \x in {0,...,1}{\draw[fill,  ] ({2*sin(180*\x+2*36)}, {2*cos(180*\x+2*36)}) circle (2pt);}
\foreach \x in {0,...,1}{\draw[fill,  ] ({2*sin(180*\x+3*36)}, {2*cos(180*\x+3*36)}) circle (2pt);}
\foreach \x in {0,...,1}{\draw[fill,  ] ({2*sin(180*\x+4*36)}, {2*cos(180*\x+4*36)}) circle (2pt);}
\draw[fill,   ] (0,1.236) circle (2pt);
\draw[fill,   ] (0,-1.236) circle (2pt);

\node at (0, -3.5){$\widetilde{G}$};

\end{tikzpicture}};

\draw[->] (T) --node[above]{$f$} (G);

\end{tikzpicture}
$$

\end{example}

We show that the classical theory of universal covers and deck transformations can be extended to homotopy covers.  
We will present these results in terms of the following groupoid notion.

\begin{definition} \label{D:gpoidcover} \cite{Higgins}  If $x$ is an object of a groupoid $A$, the star of $x$ is the set of all arrows with source $x$ and is denoted $A_x$.  Then $f:  A \to B$ is a {\bf covering of groupoids} if $f$ induces a bijection on stars $A_x \to B_{f(x)}$ for all objects $x$ of $A$.

\end{definition}

When we apply this definition to the fundamental groupoid, we get another equivalent definition of homotopy covering map.   We denote the star of $\Pi(G)$ at a vertex $v$ by $\Pi_v(G)$.  

\begin{proposition}\label{Lemma:HLiftWalk} 
 A graph homomorphism $f:  \widetilde{G} \to G$  is a homotopy covering map if and only if $f$ is a  cover as in Definition \ref{D:Cover} and the induced functor $\Pi(f):  \Pi( \widetilde{G} ) \to \Pi( G)$   is a covering of groupoids.

\end{proposition}

\begin{proof}

Suppose that $f:  \widetilde{G} \to G$ is a homotopy covering map.  To show that it induces a bijection on stars $f_*:  \Pi_{\widetilde{v}}(\widetilde{G}) \to \Pi_v(G)$, we observe that since $f$ is a covering by Observation \ref{O:useful}, Lemma \ref{L:liftwalkold} allows us to lift paths and so $f_*$ is surjective.   To see that $f_*$ is also injective, we recall that equivalences in $\Pi(G)$ are generated by prunes, anti-prunes, and  spider moves.   We again use Observation \ref{O:useful} to see that any prunable section of a walk $e \overline{e}$ lifts to a similarly prunable section $d \overline{d}$, and that any spider move that changes edges but not vertices will lift to a similar edge move  since parallel edges lift to parallel edges.  Lastly, any spider move that adjusts one vertex will lift to a  similar spider move because the bijection on second neighborhoods $N_2(v)$ respects endpoints.

Conversely, suppose $f:\widetilde{G}\to G$ is a graph cover which  induces a covering of groupoids $\Pi(f):\Pi(\widetilde{G})\to \Pi(G)$.  We know that we can lift any 2-walk uniquely by Lemma \ref{L:liftwalkold}.   This lift must respect endpoints because if two 2-walks have the same endpoint, then they represent the same arrow in $\Pi_v(G)$ and then the bijecton on stars ensures that their lifts represent the same arrow of $\Pi_{\widetilde{v}} (\widetilde{G})$ and hence have the same target vertex.

\end{proof}

\begin{definition}\label{Definition:hUC}
Given a   graph $G$ and $v\in V(G)$, we define a graph $U = U_v G$ and graph homomorphism  $\rho:  U \to G$ as follows: 

\begin{itemize}
    \item vertices of $U_v G $ are arrows in $\Pi _v G$ 
    \item  edges between arrows  $\alpha, \beta$ are defined according to the edges between their endpoints:  if $\alpha = (v v_2 \dots w)$ and $\beta = (v v_2' \dots w')$ then the edges $d$ with $\partial_0d = \alpha$ and $\partial_1d = \beta$ are given by the set   $\{d_e: e \in E(G) \textup { such that } \partial_0(e)=w, \partial_1(e)=w'\}$.
    \item the map $\rho:  U \to G$ is defined on vertices by the endpoint $\rho(\alpha) = \rho (v v_2 \dots w) =  w$ and on edges by $\rho(d_e) = e$.   
\end{itemize}
\end{definition}

\begin{example}\label{Example:HUniversalCover}
Consider the graph  $G$ from Example \ref{Example:HCover}.  Arrows in $\Pi _a G$ are equivalent if and only if they differ by spider moves. Here we use $\overline{b}$ to denote that the walk traverses through either  $b$ or $b'$, as these are  homotopic and hence identified in $\Pi _a G$.  Thus $U_aG$ may be depicted as follows:

$$\begin{tikzpicture}[scale=0.8]

\draw (-7,0)--(-5,0)--(-4,.5)--(-3,0);
\draw (-3,0) edge[bend left] (-2,0);
\draw (-3,0) edge[bend right] (-2,0);
\draw (-2,0)--(0,0)--(1,.5)--(2,0);
\draw (3,0) edge[bend left] (2,0);
\draw (3,0) edge[bend right] (2,0);
\draw (3,0)--(5,0)--(6,.5)--(7,0)--(7,0);
\draw (-5,0)--(-4,-.5)--(-3,0);
\draw (0,0)--(1,-.5)--(2,0);
\draw (5,0)--(6,-.5)--(7,0);
\node at (-7.5,0){$\cdots$};
\node at (7.5,0){$\cdots$};

\foreach \x in {-1,...,1}{ \draw[ , fill] (\x*5, 0) circle (2pt);}
\foreach \x in {-1,...,1}{ \draw[ , fill] (\x*5+1, 0.5) circle (2pt);}
\foreach \x in {-1,...,1}{ \draw[  , fill] (\x*5+1, -0.5) circle (2pt);}
\foreach \x in {-1,...,1}{ \draw[ , fill] (\x*5+2, 0) circle (2pt);}
\foreach \x in {-1,...,1}{ \draw[ , fill] (\x*5-1, 0) circle (2pt);}
\foreach \x in {-1,...,1}{ \draw[ , fill] (\x*5-2, 0) circle (2pt);}

\draw (0,0)--node[below]{\tiny $a$} (0,0);
\draw (1,0.5)--node[above]{\tiny $ab$} (1,0.5);
\draw (1,-0.5)--node[below]{\tiny $ab'$} (1,-0.5);
\draw (2,0)--node[below]{\tiny $a\overline{b}c$} (2,0);
\draw (3,0)--node[below]{\tiny $a\overline{b}cd$} (3,0);
\draw (4,0)--node[below]{\tiny $a\overline{b}cde$} (4,0);
\draw (5,0)--node[below]{\tiny $a\overline{b}cdea$} (5,0);
\draw (6,0.5)--node[above]{\tiny $a\overline{b}cdeab$} (6,0.5);
\draw (6,-0.5)--node[below]{\tiny $a\overline{b}cdeab'$} (6,-0.5);
\draw (7,0)--node[below right]{\tiny $a\overline{b}cdea\overline{b}c$} (7,0);

\draw (-1,0) --node[below]{\tiny $ae$} (-1,0);
\draw (-2,0) --node[below]{\tiny $aed$} (-2,0);
\draw (-3,0) --node[below]{\tiny $aedc$} (-3,0);
\draw (-4,0.5) --node[above]{\tiny $aedcb$} (-4,0.5);
\draw (-4,-0.5) --node[below]{\tiny $aedcb'$} (-4,-0.5);
\draw (-5,0) --node[below]{\tiny $aedc\overline{b}a$} (-5,0);
\draw (-6,0) --node[above]{\tiny $aedc\overline{b}ae$} (-6,0);
\draw (-7,0) --node[below]{\tiny $aedc\overline{b}aed$} (-7,0);

\node at (0,-1.5){$U$};

\end{tikzpicture}$$

\end{example}

\begin{obs}\label{O:parallelU}  

Recall from Notation \ref{N:unloop} that $G_s$ refers to the graph $G$ with parallel edges identified, and that Lemma \ref{L:hom}  ensures that $\Pi(G) \equiv \Pi(G_s)$.  Thus $U_vG$ has the same vertices as $U_v{G_s}$ but with additional parallel edges corresponding to those of $G$.  
\end{obs}

With this construction,  $U_v G$ is in fact a universal homotopy cover in the sense that it satisfies the appropriate universal properties.

\begin{theorem}\label{T:UMP}
    Suppose that  $f:\widetilde{G}\to G$ is a homotopy cover and   $\tilde{v}\in f^{-1}(v)$.  Then there is a unique homotopy cover  $\tilde{\rho}: U_vG \to \widetilde{G}$ such that $\rho=f\circ\tilde{\rho}$ and $\tilde{\rho}((v))=\tilde{v}$. 
\end{theorem}

\begin{proof}[Sketch of Proof]
    The proof of the above follows standard arguments as in \cite{Angluin1},\cite{Bass},\cite{Hatcher}. Any walk $\alpha$ in $G$ lifts uniquely to a walk $\widetilde{\alpha}$ starting at $\tilde{v}$ in $\widetilde{G}$, and so we define $\tilde{\rho}(\alpha)$  to be the endpoint of this lift.  This is easily extended to edges. 
\end{proof}

The previous result implies that  $U_vG, U_{\tilde{v}} \widetilde{G}$  are covers of both $G, \widetilde{G}$, so they cover each other canonically and thus are isomorphic.  Thus we also obtain the following result.  
\begin{corollary}
    If $\widetilde{G} \to G$ is a homotopy cover, then the universal homotopy cover $U$ of $G$ is also the universal homotopy cover for $\widetilde{G}$.  
\end{corollary}

\begin{example}\label{E:Ufactor}
We revisit the graphs from  Examples \ref{Example:HCover}, \ref{Example:HUniversalCover} to see how the universal homotopy covering $U \to G$ factors through the covering map $\widetilde{G} \to G$.  If we choose  $a_1\in f^{-1}(a)$, we can set  $\tilde\rho([(a)]) = a_1$.  Then the rest of $\tilde{\rho}$ follows.

$$
\begin{tikzpicture}
\node(G) at (3,0){\begin{tikzpicture}
\draw ({sin(144)},{cos(144)})\foreach \x in {3,...,6}{--({sin(72*\x)}, {cos(72*\x)})};
\draw ({sin(144)},{cos(144)}) edge[bend left] ({sin(72)},{cos(72)});
\draw ({sin(144)},{cos(144)}) edge[bend right] ({sin(72)},{cos(72)});
\draw (-.951, 0.309)--(0,0)--(.951,0.309);

\draw[ ] ({sin(72*0)}, {cos(72*0)}) -- node[above]{$b$} ({sin(72*0)}, {cos(72*0)});
\draw[  ] ({sin(72*1)}, {cos(72*1)}) -- node[right]{$c$} ({sin(72*1)}, {cos(72*1)});
\draw[ ] ({sin(72*2)}, {cos(72*2)}) -- node[below right]{$d$} ({sin(72*2)}, {cos(72*2)});
\draw[ ] ({sin(72*3)}, {cos(72*3)}) -- node[below left]{$e$} ({sin(72*3)}, {cos(72*3)});
\draw[ ] ({sin(72*4)}, {cos(72*4)}) -- node[left]{$a$} ({sin(72*4)}, {cos(72*4)});
\draw[  ] (0, 0) -- node[below]{$b'$} (0, 0);

\draw[fill,  ] ({sin(72*0)}, {cos(72*0)}) circle (2pt);
\draw[fill,   ] ({sin(72*1)}, {cos(72*1)}) circle (2pt);
\draw[fill,  ] ({sin(72*2)}, {cos(72*2)}) circle (2pt);
\draw[fill,  ] ({sin(72*3)}, {cos(72*3)}) circle (2pt);
\draw[fill,  ] ({sin(72*4)}, {cos(72*4)}) circle (2pt);
\draw[fill,   ] (0, 0) circle (2pt);

\node at (0, -1.5){$G$};

\end{tikzpicture}};

\node(T) at (-3,0){\begin{tikzpicture}[scale=0.7]
\draw ({2*sin(36*2)},{2*cos(36*2)})\foreach \x in {3,...,6}{--({2*sin(36*\x)}, {2*cos(36*\x)})};
\draw ({2*sin(7*36)},{2*cos(7*36)})\foreach \x in {8,...,11}{--({2*sin(36*\x)}, {2*cos(36*\x)})};
\draw ({2*sin(36)},{2*cos(36)}) edge[bend left]  ({2*sin(36*2)},{2*cos(36*2)});
\draw ({2*sin(36)},{2*cos(36)}) edge[bend right]  ({2*sin(36*2)},{2*cos(36*2)});
\draw ({2*sin(36*6)},{2*cos(36*6)}) edge[bend left]  ({2*sin(36*7)},{2*cos(36*7)});
\draw ({2*sin(36*6)},{2*cos(36*6)}) edge[bend right]  ({2*sin(36*7)},{2*cos(36*7)});
\draw (-1.176, 1.618)--(0,1.236)--(1.176,1.618);
\draw (-1.176, -1.618)--(0,-1.236)--(1.176,-1.618);

\draw[ ] ({2*sin(36*0)}, {2*cos(36*0)}) -- node[above]{$b_1$} ({2*sin(36*0)}, {2*cos(36*0)});
\draw[  ] ({2*sin(36*1)}, {2*cos(36*1)}) -- node[above right]{$c_1$} ({2*sin(36*1)}, {2*cos(36*1)});
\draw[ ] ({2*sin(36*2)}, {2*cos(36*2)}) -- node[right]{$d_1$} ({2*sin(36*2)}, {2*cos(36*2)});
\draw[ ] ({2*sin(36*3)}, {2*cos(36*3)}) -- node[right]{$e_1$} ({2*sin(36*3)}, {2*cos(36*3)});
\draw[ ] ({2*sin(36*4)}, {2*cos(36*4)}) -- node[below right]{$a_2$} ({2*sin(36*4)}, {2*cos(36*4)});
\draw[ ] ({2*sin(36*5)}, {2*cos(36*5)}) -- node[below]{$b_2$} ({2*sin(36*5)}, {2*cos(36*5)});
\draw[  ] ({2*sin(36*6)}, {2*cos(36*6)}) -- node[below left]{$c_2$} ({2*sin(36*6)}, {2*cos(36*6)});
\draw[ ] ({2*sin(36*7)}, {2*cos(36*7)}) -- node[left]{$d_2$} ({2*sin(36*7)}, {2*cos(36*7)});
\draw[ ] ({2*sin(36*8)}, {2*cos(36*8)}) -- node[left]{$e_2$} ({2*sin(36*8)}, {2*cos(36*8)});
\draw[ ] ({2*sin(36*9)}, {2*cos(36*9)}) -- node[above left]{$a_1$} ({2*sin(36*9)}, {2*cos(36*9)});
\draw[  ] (0,1.236) -- node[below]{$b'_1$}(0,1.236);
\draw[  ] (0,-1.236) -- node[above]{$b'_2$}(0,-1.236);

\foreach \x in {0,...,1}{\draw[fill,  ] ({2*sin(180*\x)}, {2*cos(180*\x)}) circle (2pt);}
\foreach \x in {0,...,1}{\draw[fill,   ] ({2*sin(180*\x+36)}, {2*cos(180*\x+36)}) circle (2pt);}
\foreach \x in {0,...,1}{\draw[fill,  ] ({2*sin(180*\x+2*36)}, {2*cos(180*\x+2*36)}) circle (2pt);}
\foreach \x in {0,...,1}{\draw[fill,  ] ({2*sin(180*\x+3*36)}, {2*cos(180*\x+3*36)}) circle (2pt);}
\foreach \x in {0,...,1}{\draw[fill,  ] ({2*sin(180*\x+4*36)}, {2*cos(180*\x+4*36)}) circle (2pt);}
\draw[fill,   ] (0,1.236) circle (2pt);
\draw[fill,   ] (0,-1.236) circle (2pt);

\node at (0, -3.5){$\widetilde{G}$};

\end{tikzpicture}};

\node(U) at (0,5){\begin{tikzpicture}[scale=0.8]

\draw (-7,0)--(-5,0)--(-4,.5)--(-3,0);
\draw (-3,0) edge[bend left] (-2,0);
\draw (-3,0) edge[bend right] (-2,0);
\draw (-2,0)--(0,0)--(1,.5)--(2,0);
\draw (3,0) edge[bend left] (2,0);
\draw (3,0) edge[bend right] (2,0);
\draw (3,0)--(5,0)--(6,.5)--(7,0)--(7,0);
\draw (-5,0)--(-4,-.5)--(-3,0);
\draw (0,0)--(1,-.5)--(2,0);
\draw (5,0)--(6,-.5)--(7,0);
\node at (-7.5,0){$\cdots$};
\node at (7.5,0){$\cdots$};

\foreach \x in {-1,...,1}{ \draw[ , fill] (\x*5, 0) circle (2pt);}
\foreach \x in {-1,...,1}{ \draw[ , fill] (\x*5+1, 0.5) circle (2pt);}
\foreach \x in {-1,...,1}{ \draw[  , fill] (\x*5+1, -0.5) circle (2pt);}
\foreach \x in {-1,...,1}{ \draw[ , fill] (\x*5+2, 0) circle (2pt);}
\foreach \x in {-1,...,1}{ \draw[ , fill] (\x*5-1, 0) circle (2pt);}
\foreach \x in {-1,...,1}{ \draw[ , fill] (\x*5-2, 0) circle (2pt);}

\draw (0,0)--node[below]{\tiny $a$} (0,0);
\draw (1,0.5)--node[above]{\tiny $ab$} (1,0.5);
\draw (1,-0.5)--node[below]{\tiny $ab'$} (1,-0.5);
\draw (2,0)--node[below]{\tiny $a\overline{b}c$} (2,0);
\draw (3,0)--node[below]{\tiny $a\overline{b}cd$} (3,0);
\draw (4,0)--node[below]{\tiny $a\overline{b}cde$} (4,0);
\draw (5,0)--node[below]{\tiny $a\overline{b}cdea$} (5,0);
\draw (6,0.5)--node[above]{\tiny $a\overline{b}cdeab$} (6,0.5);
\draw (6,-0.5)--node[below]{\tiny $a\overline{b}cdeab'$} (6,-0.5);
\draw (7,0)--node[below right]{\tiny $a\overline{b}cdea\overline{b}c$} (7,0);

\draw (-1,0) --node[below]{\tiny $ae$} (-1,0);
\draw (-2,0) --node[below]{\tiny $aed$} (-2,0);
\draw (-3,0) --node[below]{\tiny $aedc$} (-3,0);
\draw (-4,0.5) --node[above]{\tiny $aedcb$} (-4,0.5);
\draw (-4,-0.5) --node[below]{\tiny $aedcb'$} (-4,-0.5);
\draw (-5,0) --node[below]{\tiny $aedc\overline{b}a$} (-5,0);
\draw (-6,0) --node[above]{\tiny $aedc\overline{b}ae$} (-6,0);
\draw (-7,0) --node[below]{\tiny $aedc\overline{b}aed$} (-7,0);

\node at (0,-1.5){$U$};

\end{tikzpicture}};

\draw[->] (T) --node[below ]{$f$} (G);
\draw[->] (U) --node[above right]{$\rho$}(G);
\draw[->, dashed] (U) --node[above left]{$\exists\,\, !\,\, \tilde{\rho}$} (T);

\end{tikzpicture}
$$

\end{example}

We also have a theory of  deck transformations of homotopy covering maps for non-simple graphs.   

\begin{definition}\label{Definition:DT}
Given a graph $G$ and the universal homotopy covering map $\rho:U \to G$,  a {\bf deck transformation} is an automorphsim  $f\in \Aut(U)$ such that $\rho\circ f=\rho$.  These form a subgroup of $\Aut(U)$, which we denote by $D(G)$.   
\end{definition}
  We will show that this group is isomorphic to the isotropy group of the fundamental groupoid.  
  
\begin{definition} \label{d:fundgp}
Let $\Pi_v^vG$ denote the isotropy group of $\Pi G$ at $v$, given by the arrows of $\Pi_G$ starting and ending at $v$.  Explicitly, this is the group formed by walks that start and end at $v$ defined up to prunes and homotopy,  under concatenation.  
\end{definition}
 Because $\Pi G$ is a connected groupoid, all the isotropy groups are isomorphic, and it does not matter which vertex $v$ we choose to look at.   For $\gamma \in \Pi_v^v G$, we define a deck transformation as follows.

\begin{definition}\label{D:concact}  

    Let $\gamma\in \Pi_v^v(G)$.  Define $\psi_\gamma: U\to U$ by $\alpha\mapsto \gamma*\alpha$ for vertices $\alpha \in V(U)$, and on edges by observing that  the edges between $\alpha, \beta$ and $\gamma*\alpha, \gamma*\beta$ are both in bijection with the edge set between the endpoints $w, w'$ of $\alpha, \beta$ and so we define $\psi_\gamma(d_e) = d'_e$ for $e$ from $w$ to $w'$.   
\end{definition}

Each $\psi_\gamma$ is  deck transformation.   In fact, every deck transformation is of this form, and thus we may identify the deck transformations of $U$ with $\Pi_v^v(G)$.

\begin{theorem}\label{Theorem:DTisPi1}
The map $\Phi:  \Pi_v^v G\to D(G)$ defined by $\gamma \mapsto  \psi_\gamma$  is an isomorphism of groups.
\end{theorem}
\begin{proof}

The map  $\Phi$ is a group homomorphism because   $\Phi(\gamma*\gamma')(\alpha)=\gamma*\gamma' *\alpha  =\psi_{\gamma}\left( \psi_{\gamma'}(\alpha)\right) =\left(\Phi(\gamma)\circ \Phi(\gamma') \right)(\alpha)$.  To show that $\Phi$ is injective, suppose that  $\gamma, \gamma'\in \Pi_v^v G$ such that $\Phi(\gamma) = \Phi(\gamma')$. Then $\psi_\gamma((v)) = \psi_{\gamma'}((v))$ and so $\gamma = \gamma * (v) = \gamma'*(v) = \gamma'$. 

 To check that $\Phi$ is surjective,  we let $\varphi\in D(G)$ and define  $\gamma=\varphi((v))$.   Since $\rho(\varphi((v)))=\rho((v))= v$, the walk $\gamma$ ends at $v$ and   hence  $\Phi(\gamma)=\psi_\gamma$ is defined with  $\psi_\gamma((v))=\varphi((v))$.   
  We claim that $\varphi = \psi_\gamma$.  We prove this by induction on the length of a walk representing a vertex in $U$, starting with our base case length $0$ walk $(v)$.  Let $\alpha$ have a representative of length $n$, of the form $\beta * e$ for an edge $e$ and a length $n-1$ walk $\beta$, and assume that $\varphi(\beta) = \psi_\gamma({\beta})$.   We know that $e$ defines an edge $d_e$ from $\beta$ to $\alpha$ in $U$, and since $\varphi$ is a graph homomorphism that commutes with $\rho$ we have an edge $\varphi(d_e) = d'_e$ from $\varphi(\beta) $ to $\varphi(\alpha)$.  Therefore $\varphi(\alpha) = \varphi(\beta) * d'_e = \gamma * \beta * d'_e $.   Uniqueness of path lifting ensures that $\gamma * \beta * d'_e = \gamma * \alpha$.   Thus   $\varphi = \psi_\gamma$ and so $\Phi$ is surjective, completing the proof.

\end{proof}

   Given a subgroup $S \leq D(G)$, we define $\widetilde{G} = U/S$ as the quotient graph where we identify vertices $[\alpha] = [s\alpha]$ and edges $[sd_e ] = [d_e]$  for any  $s \in S$.  We can show that any cover is the quotient of $U$ by some subgroup of the deck transformation group $D(G)$.  Our proof uses  groupoid coverings and Proposition \ref{Lemma:HLiftWalk}.  We start by showing that the resulting quotient is a homotopy cover.

\begin{theorem}\label{Theorem:QuotientisCover}
Let $G$ be a graph with universal homotopy cover $ \rho: U \to G$, and let $S\leq D(G)$.   Define $r:U/S\to G$ on vertices by  $r[\alpha] = \rho(\alpha)$ and on edges $r[d_e] = \rho(d_e)=e$.    Then $r$ is a homotopy covering map.  

\end{theorem}

\begin{proof}
We first note that $r$ is well-defined on equivalence classes forming the vertices  $U/S$ because  for any $s \in S\leq D(G)$, we have   $\rho s \alpha = \rho \alpha$;  and hence $r$ is surjective on vertices because $\rho$ is.   The map $r$ is similarly well-defined on edges, and moreover   since $\rho (s d_e)  =  \rho(d_e) = e$ the edges $[d_e]$ from $[\alpha]$ to $[\beta]$  are still in bijection with edges $e$ from their endpoints $w, w'$  in $G$.   Thus $r$ induces a bijection from $N_E([\alpha]) $ to  $N_E(w)$.    Therefore $r$ is a cover as in Definition \ref{D:Cover}.

We now show that $r$ induces a covering map on fundamental groupoids $\Pi(r):  \Pi(U/S) \to \Pi(G)$, meaning that it is bijective on stars. 
By the definition of the map $r$, $r p = \rho$ where $p$ is the projection map $p:  U \to U/S$.   Functoriality ensures that $\Pi(\rho ) = \Pi(rp) = \Pi(r) \Pi(p)$.   Now  $\Pi(p)$ is surjective on stars:  given a walk $([d_1][d_2]\cdots [d_n])$ in $U/S$, we know that there are representatives $d_i$ in $U$ and $s_i \in S$ such that  $\partial_1(d_{i-1}) =\partial_0(s_id_i)$.  By   applying the $s_i$ in turn to the end of the walk,  we obtain a walk in  $U$  that maps to $([d_1][d_2]\cdots [d_n])$ by $p$.       

Since $\Pi(p)$ is surjective on stars and $\Pi(\rho) = \Pi(r) \Pi(p)$ is bijective on stars, we conclude that  $\Pi(r)$ is bijective on stars. 
By Proposition \ref{Lemma:HLiftWalk}, $r$ is a homotopy covering map.

\end{proof}

\begin{proposition}\label{Proposition:Coversmodthemselves}
Let $G$ be a graph, and $U$ be its universal cover.  Then  any connected homotopy cover $\widetilde{G}$  is of the form  $U/S \cong \widetilde{G}$ for some subgroup $S \leq D(G)$. 
\end{proposition}
\begin{proof}
Suppose that $f:  \widetilde{G} \to G$ is a homotopy covering map.   
Theorem \ref{T:UMP} states that $U$ is also the universal homotopy cover of $\widetilde{G}$ and there exists a homotopy covering map $\tilde{\rho}:  U \to \widetilde{G}$ such that $\rho = pf$.  Define   $S = \{s \in D(G) : s\tilde{\rho} = \tilde{\rho}\}$. We wish to show that  $U/S \cong \widetilde{G}$.   In fact, $S = D(\widetilde{G})$, and  thus it suffices to  consider the case  $S=D(G)$ and prove that $U/D(G)\cong G$.

Let $r:U/D(G)\to G$ be defined by $r[u]=\rho(u)$.  Then $r$ is a homotopy covering map by Theorem \ref{Theorem:QuotientisCover}, and hence surjective.    We will show that  it is also injective. Suppose that  $r[\alpha]=w=r[\beta]$ for vertices $\alpha, \beta \in U$.   Since $\alpha, \beta$ have the same endpoint, define $\gamma = \alpha * \beta^{-1}$, giving a walk starting and ending at $v$, and let $s  = \psi_\gamma$ as in Definition \ref{D:concact}.   Then $\psi_\gamma \in D(G)$ and $\psi_\gamma(\beta) = \alpha$ and so  $[\alpha] = [\beta]$ in $U/D(G)$. Therefore $r$ is a bijection of vertices, and also a cover which gives a bijection on neighborhoods.   So $r$ is also a bijection on edges and hence an isomorphism.  

\end{proof}

\begin{example}\label{E:Quotient} We consider the graphs $G, \widetilde{G}, U$ and the corresponding maps from Example \ref{E:Ufactor}.  One may  verify that since $G\simeq C_5$, it follows that $\Pi_a^aG\cong \mathbb{Z}$, generated by $(abcdea)$.  Moreover,  $\tilde{G}\cong U/\langle (abcdea)^2 \rangle$ in accordance with Proposition \ref{Proposition:Coversmodthemselves}.

\end{example}

\begin{obs}\label{O:parallelcover}
Since all homotopy covers of $G$ are isomorphic to quotients of $U$,  Observation \ref{O:parallelU} ensures that any  homotopy cover $\widetilde{G}$ of a loop-free graph $G$ is isomorphic to a supergraph of a homotopy cover $\widetilde{G_s}$ of $G_s$, where vertices $[\alpha], [\beta]$ in $\widetilde{G}$ have the same number of parallel edges between them as $r[\alpha], r[\beta]$ do.  
\end{obs}

We can use  Proposition \ref{P:reflexivegroupoid} to give a description of the universal homotopy cover of a reflexive graph as it relates to the universal homotopy cover of its unlooped subgraph.

\begin{notation} \label{N:2cover}
    Let $G$ be a  reflexive graph.  Let  $f:C\to G_u$ be a homotopy covering map of the unlooped graph  $G_u$.  We define $C_\ell$ to be the graph $C$ with loops appended so that for any vertex, $\alpha$, we have an isomorphism between the loops on $\alpha$ and the loops on the vertex $f(\alpha)$ in the original graph $G$.  
\end{notation}

With this definition, we get the following.  
\begin{theorem}\label{T:reflexiveU}

Suppose that $G$ is a reflexive graph.  We define the cover $C$ of its unlooped graph $G_u$ by $C = U_vG_u/T$, the quotient of the universal cover of $G_u$ by the the subgroup $T \leq D(G_u) \cong \Pi_v^v(G_u)$ generated by 3-cycles.   Then the universal cover of $G$ satisfies $$U_vG \cong C_\ell \times K_2.$$  
\end{theorem}
 \begin{proof}
We know that $U_v G$ has vertices which are defined by walks from $v$ in $\Pi_v(G)$.  We apply the isomorphism of Proposition \ref{P:reflexivegroupoid} to get an arrow in the groupoid $\Pi_v(G_u)/T \times \zed/2$  starting at $v$, giving an isomorphism between the vertices of $U_v G$ and $C_\ell\times K_2$.  We then define the mapping $U_vG \to C_\ell\times K_2$ on edges.  Any edge $d_e\in(U_v G)$ corresponds to extending a walk $\alpha$ by a single edge $e\in E(G)$ to $\alpha * e$.   In this case, the parity of $\alpha$ and $\alpha * e$ are opposite, and so $\alpha$ and $\alpha*e$ are mapped to $(\alpha_u, i)$ and $((\alpha*e)_u, 1-i)$ in $C_\ell\times K_2$.    If $\alpha_u = (\alpha * e)_u$ then $e = \ell $ was a loop, and we define the image of $d_e$ to be the edge   $c_e= d_\ell$ with $\partial_0(c_\ell)=(\alpha_u, i) , \partial_1(c_\ell)=(\alpha_u, 1-i)$.   Otherwise we have $\alpha_u * e = (\alpha * e)_u$ and there is an edge $c_e$,  $\partial_0(c_e)=(\alpha_u , i), \partial_1(c_e)=  (\alpha_u * e, 1-i)$. We then map $d_e\mapsto c_e$.  It is straightforward to check that this gives a bijective map on edges in which $\overline{d_e}\mapsto \overline{c_e}$.  Thus we get an isomorphism of graphs.

 \end{proof}

\begin{corollary}
    Let $G$  be reflexive and $C_3$-free graph, and let $U_v G$ denote its universal homotopy cover. Then $U_v G\cong U_v G_u\square K_2$.
\end{corollary}
\begin{proof}  We know that $U_v G \cong C_\ell\times K_2$ and 
     since  $G$ is $C_3$ free, $U_v G_u=C$ and so the parity of a walk in $C$ is well-defined.    Define a homomorphism  $\Phi:  C_\ell\times K_2 \to C\square K_2$ by $(\gamma, i) \to (\gamma, i+k)$ where $k$ is the length of $\gamma$ (mod 2).   Since $\gamma$ and $\gamma*e$ always have opposite parity, the vertices $(\gamma, i)$ and $(\gamma*e, 1-i)$ are always mapped to the same layer in the box product and this gives a graph homomorphism.      It is straightforward to see that $\Phi$ is bijective on vertices and edges and thus is an isomorphism.

\end{proof}

\begin{example}\label{E:Ucover}
We consider the graph  $G$ from Example \ref{E:reflexivegroupoid}. In this case, the cover $C = U G_u/T$ is a quotient of the universal homotopy cover by the triangle $(axea)$.   The result is drawn below,  with covering map $f$ defined by $v_i \to v$:  
\[
\begin{tikzpicture}[scale=0.8]
    \draw (-2,0)--(0,0);
    \draw (10,0)--(12,0);
    \foreach \x in {0,..., 1}{ 
    \draw ({\x*6},0)--({\x*6+1},0);
  \draw ({\x*6+2},0)--({\x*6+6},0);
  \draw ({\x*6+1},0) edge[bend left] ({\x*6+2},0);
  \draw ({\x*6+1},0) edge[bend right] ({\x*6+2},0);
    \draw[fill, ] ({\x*6+1}, 0) node[below]{$b_{\x}$}  circle (2pt);
\draw[fill, ] ({\x*6}, 0) node[below]{$a_{\x}$}  circle (2pt);
\draw[fill,   ] ({\x*6+2}, 0) node[below]{$c_{\x}$}  circle (2pt);
\draw[fill, ] ({\x*6+3}, 0) node[below]{$d_{\x}$}  circle (2pt);
\draw[fill,  ] ({\x*6+4}, 0) node[below]{$e_{\x}$}  circle (2pt);
\draw[fill] ({\x*6+5}, 0) node[below]{$x_{\x}$}  circle (2pt);

}

\draw[fill,  ] ({-2}, 0) node[below]{$e_{-1}$}  circle (2pt);
\draw[fill] ({-1}, 0) node[below]{$x_{-1}$}  circle (2pt);
\draw[fill, ] ({12}, 0) node[below]{$a_{2}$}  circle (2pt);

\draw (-2,0) edge[bend left] (0,0);
\draw (4,0) edge[bend left] (6,0);
\draw (10,0) edge[bend left] (12,0);

\node at (-3,0){$\cdots$};
\node at (13,0){$\cdots$};

\node at (5,-1){$UG_u / T$};

\end{tikzpicture}
\]  
Following Theorem \ref{T:reflexiveU}, we obtain $U G= C_\ell\times K_2$ with vertices given by  ordered  pairs $(\alpha, i)$ where $\alpha\in V(C), i\in \mathbb{Z}/2$ and edges defined by  $(\alpha, 0)\sim (\beta, 1)$ for each edge from $f(\alpha)$ to $f(\beta)$ in $G$.

\[\begin{tikzpicture}

\foreach \x in {-1,..., 10}{ 
    \ifthenelse{\NOT -6 = \x \AND \NOT 0 = \x \AND \NOT 6 = \x}{\draw(\x,0) -- (\x+1,2); \draw(\x,2) -- (\x+1,0);}{
    \draw(\x,0) edge[bend left] (\x+1,2);
    \draw(\x,0) edge[bend right] (\x+1,2);
    \draw(\x,2) edge[bend left] (\x+1,0);
    \draw(\x,2) edge[bend right] (\x+1,0);
    }
}
\foreach \x in {-1,..., 11}{ 
\draw(\x,0) -- (\x,2);
}

\foreach \x in {1,..., 0}{ \draw[fill, ] ({\x*6}, 0) node[below]{\tiny$(b_{\x},0)$}  circle (2pt);
\draw[fill, ] ({\x*6-1}, 0) node[below]{\tiny$(a_{\x},0)$}  circle (2pt);
\draw[fill,   ] ({\x*6+1}, 0) node[below]{\tiny$(c_{\x},0)$}  circle (2pt);
\draw[fill, ] ({\x*6+2}, 0) node[below]{\tiny$(d_{\x},0)$}  circle (2pt);
\draw[fill,  ] ({\x*6+3}, 0) node[below]{\tiny$(e_{\x},0)$}  circle (2pt);
\draw[fill] ({\x*6+4}, 0) node[below]{\tiny$(x_{\x},0)$}  circle (2pt);

\draw[fill, ] ({\x*6-1}, 2) node[above]{\tiny$(a_{\x},1)$}  circle (2pt);
\draw[fill, ] ({\x*6}, 2) node[above]{\tiny$(b_{\x},1)$}  circle (2pt);
\draw[fill,   ] ({\x*6+1}, 2) node[above]{\tiny$(c_{\x},1)$}  circle (2pt);
\draw[fill, ] ({\x*6+2}, 2) node[above]{\tiny$(d_{\x},1)$}  circle (2pt);
\draw[fill,  ] ({\x*6+3}, 2) node[above]{\tiny$(e_{\x},1)$}  circle (2pt);
\draw[fill] ({\x*6+4}, 2) node[above]{\tiny$(x_{\x},1)$}  circle (2pt);

\draw ({\x*6+3}, 2)  -- ({\x*6+5}, 0);
\draw ({\x*6+3}, 0)  -- ({\x*6+5}, 2);
\draw ({\x*6+3}, 2)  -- ({\x*6+4}, 0);
\draw ({\x*6+3}, 0)  -- ({\x*6+4}, 2);

}

\draw[fill, ] ({11}, 0) node[below]{\tiny$(a_2,0)$}  circle (2pt);

\draw[fill, ] (11, 2) node[above]{\tiny$(a_2,1)$}  circle (2pt);

\foreach \x in {0,..., 2}{
\draw (12,\x) -- node{$\cdots$} (12,\x);
\draw (-2,\x) -- node{$\cdots$} (-2,\x);
}

\node at (6, -1){$U G$};

\end{tikzpicture}\]

\end{example}

\begin{proposition} Let $G$ be a reflexive graph.   Then any homotopy cover of $G$ is related to a homotopy cover $C$ of the unlooped graph $G_u$ that lifts 3-cycles to 3-cycles, and has one of the following forms: \begin{itemize}
    \item the reflexive graph $C_\ell$,    or
    \item  an unlooped graph of the form $C_\ell \times K_2$.
\end{itemize}

\end{proposition}
\begin{proof} 
By Proposition \ref{Proposition:Coversmodthemselves}, any homotopy cover $\widetilde{G}$ of $G$ has the form $\widetilde{G}\cong U_vG/S$ for $S \leq D(G)$.  And by Theorem \ref{Theorem:DTisPi1} and Proposition \ref{P:reflexivegroupoid},   $D(G)\cong \Pi_v^v G\cong \Pi_v^v G_u/T\times \mathbb{Z}/2$.  So $S$ has the form $S' \times \zed/2 $ or $S' \times \{0\}$, where $S'$ is a subgroup of $\Pi_v^vG_u$ which contains $T$.    We can define the homotopy cover of the unlooped graph to be $C = UG_u / S'$.  Then  if $S$ is of the form $S' \times \zed/2 $ we obtain a reflexive homotopy cover $\widetilde{G} $ of the first form, and if  $S$ is of the form $S' \times \{0\}$ we obtain $\widetilde{G} $ containing no loops of the second form.

\end{proof}

\begin{example}
    Recall $G$ and its fundamental groupoid and universal cover from Examples \ref{E:reflexivegroupoid} and \ref{E:Ucover}. Below, we depict  two 2-covers of $G$, one reflexive and one not.
    \[\begin{tikzpicture}
        \draw ({2*sin(2*30)}, {2*cos(2*30)})\foreach \x in {3,...,7}{--({2*sin(\x*30)}, {2*cos(\x*30)})};
        \draw ({2*sin(8*30)}, {2*cos(8*30)})\foreach \x in {9,...,13}{--({2*sin(\x*30)}, {2*cos(\x*30)})};
        \foreach \x in {0,1}{
            \draw[fill, ] ({2*sin(\x*180)}, {2*cos(\x*180)}) node[below]{$a_{\x}$}  circle (2pt);
            \draw[fill, ] ({2*sin(\x*180+30)}, {2*cos(\x*180+30)}) node[below]{$b_{\x}$}  circle (2pt);
            \draw[fill, ] ({2*sin(\x*180+60)}, {2*cos(\x*180+60)}) node[below]{$c_{\x}$}  circle (2pt);
            \draw[fill, ] ({2*sin(\x*180+90)}, {2*cos(\x*180+90)}) node[below]{$d_{\x}$}  circle (2pt);
            \draw[fill, ] ({2*sin(\x*180+120)}, {2*cos(\x*180+120)}) node[below]{$e_{\x}$}  circle (2pt);
            \draw[fill, ] ({2*sin(\x*180+150)}, {2*cos(\x*180+150)}) node[below]{$x_{\x}$}  circle (2pt);
            \draw ({2*sin(\x*180)}, {2*cos(\x*180)}) edge[bend left]({2*sin(\x*180-60)}, {2*cos(\x*180-60)}) ;
            \draw ({2*sin(\x*180+30)}, {2*cos(\x*180+30)}) edge[bend left] ({2*sin(\x*180+60)}, {2*cos(\x*180+60)});
            \draw ({2*sin(\x*180+30)}, {2*cos(\x*180+30)}) edge[bend right] ({2*sin(\x*180+60)}, {2*cos(\x*180+60)});

        }
        \foreach \x in {0,...,11}{
            \draw ({2*sin(\x*30)}, {2*cos(\x*30)}) to [in=50,out=140,loop] ({2*sin(\x*30)}, {2*cos(\x*30)});
        }
    \end{tikzpicture}
\ \ \ \ \ \ \ \ \ \ \ \ 
    \begin{tikzpicture}
        \foreach \x in {0,...,5}{
            \draw ({2*sin(\x*60)}, {2*cos(\x*60)}) -- ({sin(\x*60)}, {cos(\x*60)});
            
        }
        \foreach \x in {2,...,6}{
            \draw ({2*sin(\x*60)}, {2*cos(\x*60)}) -- ({sin(\x*60+60)}, {cos(\x*60+60)});
        }
        \foreach \x in {2,...,6}{
            \draw ({2*sin(\x*60+60)}, {2*cos(\x*60+60)}) -- ({sin(\x*60)}, {cos(\x*60)});
        }
        \foreach \x in {0,...,1}{ 
            \draw[] ({(\x+1)*sin(60)}, {(\x+1)*cos(60)}) edge[bend right] ({(2-\x)*sin(120)}, {(2-\x)*cos(120)});
            \draw[] ({(\x+1)*sin(60)}, {(\x+1)*cos(60)}) edge[bend left] ({(2-\x)*sin(120)}, {(2-\x)*cos(120)});
            \draw[] ({(\x+1)*sin(60)}, {(\x+1)*cos(60)}) edge[bend right] ({(2-\x)*sin(120)}, {(2-\x)*cos(120)});

            \draw[fill, ] ({(\x+1)*sin(0)}, {(\x+1)*cos(0)}) node[below]{$a_{\x}$}  circle (2pt);
            \draw[fill, ] ({(\x+1)*sin(60)}, {(\x+1)*cos(60)}) node[below]{$b_{\x}$}  circle (2pt);
            \draw[fill, ] ({(\x+1)*sin(120)}, {(\x+1)*cos(120)}) node[below]{$c_{\x}$}  circle (2pt);
            \draw[fill, ] ({(\x+1)*sin(180)}, {(\x+1)*cos(180)}) node[below]{$d_{\x}$}  circle (2pt);
            \draw[fill, ] ({(\x+1)*sin(240)}, {(\x+1)*cos(240)}) node[below]{$e_{\x}$}  circle (2pt);
            \draw[fill, ] ({(\x+1)*sin(300)}, {(\x+1)*cos(300)}) node[below]{$x_{\x}$}  circle (2pt);
        }
        \draw[] ({1*sin(0)}, {1*cos(0)}) edge[bend right=70] ({2*sin(240)}, {2*cos(240)});
        \draw[] ({2*sin(0)}, {2*cos(0)}) edge[bend left=70] ({1*sin(240)}, {1*cos(240)});
    \end{tikzpicture}\]
    To the left, we have $U/\langle (abcdexa)^2\rangle \times \zed/2$, and to the right $U/\langle (abcdexa) \rangle \times \{0 \}$.
\end{example}

We finish out our discussion of homotopy covers by showing that the they satisfy a general homotpy lifting property.  

\begin{theorem}\label{T:HLiftGeneral}
Let  $ f: \widetilde{G}\to G$ be a homotopy covering map, and suppose we have  $\phi, \psi:K\to G$ which are homotopic via a homotopy $H:K \times I_n \to G$.    Suppose that we can lift $\phi$ to a map  $\tilde{\phi}:K\to \widetilde{G}$ such that $\phi=f \circ \widetilde{\phi}$.  Then there is a lift $\widetilde{H}:K \times I_n \to \widetilde{G}$ such that ${\tilde{H}}|_{K \times \{0\}} = \widetilde{\phi}$   and $ f \circ \widetilde{H} = H$.  In particular, this will ensure that ${\tilde{H}}|_{K \times \{n\}} = \widetilde{\psi}$ is a lift of $\psi$ such that $\widetilde{\phi} \simeq \widetilde{\psi}$.  
\[ \xymatrix{ && \widetilde{G}  \ar[dd]^{f}  \\
 K \times \{ 0 \} \ar@{^{(}->}[d]\ar[urr]^{\widetilde{\phi}} & & \\
 K \times I_n \ar[rr]^{H:  \phi \simeq \psi} \ar@{-->}[uurr]_{\widetilde{H}} && G  } \]
\end{theorem}

\begin{proof}
We create the lift one step at a time through the length of the homotopy, so it is sufficient to consider the case when  $\phi$ and $\psi$ are connected by a length $1$ homotopy.  So consider $K \times I_1$.  We have $\widetilde{H}$ defined on $K \times \{ 0\}$.   We first extend it to the edges in $K \times I_1$ of the form $(k, e_1)$ where $e_1$ is the edge of $I_1$ with $\partial_0(e_1) = 0$ and $\partial_1(e_1) = 1$.   We do this by lifting the 2-walk $\phi(k, e_1) \phi(\overline{k}, \overline{e}_1)$ starting at the vertex $\widetilde{\phi}(\partial_0(k), 0)$.   By our observation, this lift is a walk in $\widetilde{G}$ of the form  $(d \overline{d})$.   We define $\widetilde{H}(k, e_1) = d$.   To ensure our map respects involution we define $\widetilde{H}(\overline{k}, \overline{e}_1) =  \overline{d}$. 

This means that we must define $\widetilde{H}(\partial_1k, 1) = \partial_1 d$.   This is well-defined, since if a vertex $v =\partial_1k = \partial_1k'$ then the 2-walks $H(\overline{k}, \ell_0) H(k, e_1)$ and $H(\overline{k}', \ell_0) H(k', e_1)$ lift to $\widetilde{\phi}(\overline{k}, \ell_0) d$ and $\widetilde{\phi}(\overline{k}', \ell_0) d'$ and so $\partial_1 d   = \partial_1 d'$.  

Lastly we define $\widetilde{H} $ on the edges $(k, \ell_1)$ by lifting the 2-walks $H(\overline{k}, e_1) H(k, \ell_1) $ to $(\overline{c} d)$ and setting $\widetilde{H}(k, \ell_1) = d$.    Observe that if $\partial_0 k = v$ and $\partial_1k = w$  then  $\partial_0 d = \widetilde{H} (v, 1)$  since $\partial_1 \overline{c} = \widetilde{H}(v, 1)$ by definition, and $\partial_1 d= \widetilde{H} (w, 1)$ since the lifted 2-walk must end in the same place as the lifted 2-walk $H(\overline{k}, \ell_0) H(k, e_1)$.   Lastly, $\widetilde{H}(\overline{k}, \ell_1)   = \overline{d} = \overline{\widetilde{H}(k, \ell_1)}$ by uniqueness of lifts.  

\end{proof}

\section{Future Directions}

In \cite{Docht2}, Dochtermann  presents a fundamental group for the category of reflexive graphs, consisting of classes of closed looped walks under $\times$-homotopy.  This fundamental group, although also defined via $\times$-homotopy, differs from the one presented in \cite{CS2} and this paper.  This group can readily be extended to a groupoid, and it should  be possible to develop a theory of universal covers and deck transformations in the category of reflexive graphs.  An open question is how such a construction would relate to the universal cover of this paper given for reflexive graphs as part of the more general graph category.   

In \cite{Angluin1}, Angluin  posed a conjecture that any two finite graphs which shared a universal cover would also share a common finite cover.  She proved some partial results \cite{Angluin2}, and the conjecture was finally resolved by Leighton \cite{Leighton}, and extended in a number of ways \cite{bridsonshepherd}, \cite{shepherd1}, \cite{woodhouse2},  \cite{woodhouse2018revisiting} .  The establishment of a universal $\times$-homotopy cover begs the same question: if two finite graphs shared a universal homotopy cover, would they necessarily share a  common finite homotopy cover?

This paper generalizes the notions of $\times$-homotopy and coverings which lift $\times$-homotopy to the most general graph category: that which allows but doesn't require loops and multiple parallel edges.  The next level of generalization would be to extend these ideas to hypergraphs, using the same interval graph and the product in the category of hypergraphs to establish an analogous  definition to Definition \ref{D:htpy}.  It would be interesting to investigate how the notions of $\times$-homotopy, fundamental groupoids, and $\times$-homotopy lifting covers behave in this setting, and what is gained or lost when restricted to graphs.

\section*{Acknowledgement} 
We wish to thank the referee for their helpful comments.

\end{document}